\documentclass[12pt,a4paper]{article}
\usepackage{amsmath,amssymb,amscd,amsthm}

\theoremstyle{plain}
\newtheorem{theo}{Theorem}
\newtheorem{prop}[theo]{Proposition}
\newtheorem{coro}[theo]{Corollary}
\newtheorem{lem}[theo]{Lemma}
\newtheorem{conj}{Conjecture}
\newtheorem{que}[conj]{Question}

\theoremstyle{definition}
\newtheorem{defi1}{Definition}
\newenvironment{defi}[1][]{\begin{defi1}[#1]\rm}{\end{defi1}}

\newtheorem{exe}{Example}

\newenvironment{ex}[1][]{\begin{exe}\rm}{\end{exe}}

\newtheorem{notati}[defi1]{Notation}

\newtheorem{rema}{Remark}

\newcommand{\pf}{ {\noindent {\bf Proof}}\ }

\newcommand{\mythmname}{}
\newtheoremstyle{mytheorem}% name
	{3pt}% Space above
	{3pt}% Space below
	{\it}% Body font
	{}% Indent amount
	{}% Theorem head font
	{{\bf .}}% Punctuation after theorem head
	{.5em}% Space after theorem head
	{\mythmname }% 
\theoremstyle{mytheorem}
\newtheorem{namedtheorem}{Name}
\newcommand{\renametheorem}[1]
{
	\renewcommand{\mythmname}{{\bf #1}}
}

\newcommand{\A}{\ensuremath{\mathbb A}}

\newcommand{\C}{\ensuremath{\mathbb C}}

\newcommand{\K}{\ensuremath{\mathbb K}}
\newcommand{\Q}{\ensuremath{\mathbb Q}}

\newcommand{\s}{\ensuremath{\sigma}}

\newcommand{\pri}{\ensuremath{\smallsetminus}}

\DeclareMathOperator{\Spec}{Spec}

\begin{document}

\begin{center}
 {\LARGE\bf Some Families of Polynomial Automorphisms~III}\\
\vspace{.4cm}
{\large Eric Edo\footnotemark, Drew Lewis\footnotemark}
\end{center}

{\small {\bf Abstract.}  We prove that the closure (for the Zariski topology)
of the set of polynomial automorphisms of the complex affine plane whose polydegree is $(cd-1,b,a)$
contains all triangular automorphisms of degree $cd+a$ where $a,b\ge 2$ and $c\ge 1$ are integers and $d=ab-1$.
When $b=2$, this result gives a family of counterexamples to a conjecture of Furter.

\footnotetext[1]{ERIM, University of New Caledonia; eric.edo@univ-nc.nc}
\footnotetext[2]{Department of Mathematics, University of Alabama; dlewis@ua.edu}

\section{Introduction}

\hspace*{.6cm} Let $\K$ be a field. We denote by ${\cal G}(\K)$ the group of polynomial automorphisms of the affine plane
$\A^2_\K={\rm Spec}(\K[X,Y])$. An element $\s\in {\cal G}(\K)$ is defined by
a pair of polynomials $(f,g)\in\K[X,Y]^2$ such that
$\K[f,g]=\K[X,Y]$, and we set $\s=(f,g)$. We define the degree of
$\s\in {\cal G}(\K)$ by $\deg(\s)=\max\{\deg(f),\deg(g)\}$. We denote by ${\cal A}(\K)$
the subgroup of \textit{affine} automorphisms (\emph{i.e.} automorphisms of degree $1$) and
by ${\cal B}(\K)$  the subgroup of \textit{triangular} automorphisms (of
the form $(aX+P(Y),bY+c)$ with $a,b\in\K^*$, $c\in\K$ and $P\in\K[Y]$).
For a general reference on polynomial automorphisms, see~\cite{vdE} or \cite{W}.

The classical Jung-van der Kulk theorem (\cite{J} and~\cite{vdK}, or~\cite{vdE} for a modern treatment) gives ${\cal G}(\K)$ the structure of an amalgamated free product of ${\cal A}(\K)$ and ${\cal B}(\K)$ along ${\cal A}(\K)\cap{\cal B}(\K)$.
This property allows us to define the {\it polydegree} of $\s\in{\cal G}(\K)$ as the (unique) sequence
of the degrees of the triangular automorphisms
in a decomposition of $\s$ as a  product of affine and triangular
automorphisms (cf.~\cite{FM}). We denote by ${\cal G}(\K)_d$ the set of all
automorphisms of ${\cal G}(\K)$ whose polydegree is $d=(d_1,\ldots,d_l)$
where $d_1,\ldots,d_l\ge 2$ are integers. We say the length of the sequence $(d_1,\ldots,d_l)$ is $l$.
As a convention, we set ${\cal G}(\K)_{\emptyset}={\cal A}(\K)$
where $\emptyset$ is the empty sequence. If ${\cal D}$ denotes the set of sequences of integers $\ge 2$,
including the empty sequence,
$${\cal G}(\K)=\coprod_{d\in{\cal D}}{\cal G}(\K)_d.$$
We use those notations mostly in the case $\K=\C$ and we denote simply ${\cal G}={\cal G}(\C)$, ${\cal A}={\cal A}(\C)$,
${\cal B}={\cal B}(\C)$ and ${\cal G}_d={\cal G}(\C)_d$.

The group ${\cal G}$ can be endowed with the structure of an
infinite-dimensional algebraic variety (cf.~\cite{S}).
If ${\cal H}\subset {\cal G}$, we denote by $\overline{{\cal H}}$ the
closure of ${\cal H}$ in ${\cal G}$ for the Zariski topology associated with
this structure. To check if an automorphism $f\in{\cal G}$ is in $\overline{{\cal G}_d}$, we use the Valuation Criterion
due to Furter (see Corollary 1.1 in~\cite{F3}):

\begin{theo}[Valuation Criterion]\label{thm:valcri}
Let $\tau\in{\cal G}$ be an automorphism and let $d\in{\cal D}$ be a polydegree.
Then $\tau\in\overline{{\cal G}_d}$ if and only if $\tau=\lim_{Z\to 0}\s_Z$ for some $\s_Z\in{\cal G}_d(\C((Z)))$.
\end{theo}

In the previous theorem $\C((Z))$ is the field of fractions of $\C[[Z]]$, the ring of formal power series in $Z$.
In this paper, we use the Valuation Criterion only to prove that some $\tau\in{\cal B}$ belongs to $\overline{{\cal G}_d}$
and the $\s_Z$ appearing in our proofs are always elements of ${\cal G}_d(\C[Z])$. In this context,
$\lim_{Z\to 0}\s_Z$ is simply the image of $\s_Z$ modulo $Z$.

It is natural to to examine the interaction between the structure of $\mathcal{G}$ as an infinite dimensional variety and the amalgamated free product structure of $\mathcal{G}$.  The general question we are interested in is (see~\cite{F3}):
\begin{que}
Let $d=(d_1,\ldots,d_l)$ and $e=(e_1,\ldots,e_m)$ be two degree sequences in ${\cal D}$.
What conditions on $d$ and $e$ guarantee that $\mathcal{G}_d\subset\overline{\mathcal{G}_e}$?
\end{que}

It's clear that $\mathcal{G}_d\subset\overline{\mathcal{G}_e}$ is equivalent to $\overline{\mathcal{G}_d}\subset\overline{\mathcal{G}_e}$
and thus this relation induces an order on ${\cal D}$. Following notations from~\cite{F3}, we sometimes denote this order by $d\sqsubseteq e$ .

An obvious necessary condition for the inclusion $\mathcal{G}_d\subset\overline{\mathcal{G}_e}$ is 
$\mathcal{G}_d\cap\overline{\mathcal{G}_e}\ne\emptyset$. 
In the case $l=m$, this condition is also sufficient (Theorem~C in~\cite{F3}):

\begin{theo}[Furter]
Let $d=(d_1,\ldots,d_l)$ and $e=(e_1,\ldots,e_l)$ be two degree sequences in ${\cal D}$ with the same length.
The following assertions are equivalent:\\
(i) $\mathcal{G}_d\subset\overline{\mathcal{G}_e}$\hspace{1cm}
(ii) $\mathcal{G}_d\cap\overline{\mathcal{G}_e}\ne\emptyset$\hspace{1cm}
(iii) $d_i\le e_i$ for all $i\in\{1,\ldots,l\}$.
\end{theo}

Friedland and Milnor (\cite{FM}) proved that $\mathcal{G}_d$ is an analytic variety of dimension $d_1+\cdots+d_l+6$.
In the case $d\ne e$, we have $\mathcal{G}_d\cap\mathcal{G}_e=\emptyset$ and the dimension constraint implies  
that $d_1+\cdots+d_l \leq e_1+\ldots+e_m-1$ is a necessary condition for the inclusion $\mathcal{G}_d\subset\overline{\mathcal{G}_e}$
which is equivalent to $\mathcal{G}_d\subset\overline{\mathcal{G}_e}\pri \mathcal{G}_e$.
We can ask a more precise question:

\begin{que}\label{ques:dim}
Let $d=(d_1,\ldots,d_l)$ and $e=(e_1,\ldots,e_m)$ be two degree sequences in ${\cal D}$
such that $\mathcal{G}_d\cap\overline{\mathcal{G}_e}\ne\emptyset$ and $d_1+\cdots+d_l \leq e_1+\ldots+e_m-1$.
Then is $\mathcal{G}_d\subset\overline{\mathcal{G}_e}$?
\end{que}

One necessary condition is due to Furter (\cite{F2} Theorem~1), who showed that the length of an automorphism is lower semicontinuous.  This result can be reformulated in the following way:

\begin{theo}[Furter]
Let $d=(d_1,\ldots,d_l)$ and $e=(e_1,\ldots,e_m)$ be two degree sequences in ${\cal D}$.
If ${\cal G}_d\cap \overline{{\cal G}_e}\ne\emptyset$ then $l\le m$. 
\end{theo}

In the case $l=1$ and $m=2$, we usually (see~\cite{F2}) make the following conjecture (which implies that
the answer to Question~\ref{ques:dim} is ``yes''):
\begin{conj}\label{conj:dim2}
Let $a,b,c\ge 2$ be integers.\\
A) If $c\le a+b-1$ then ${\cal G}_{(c)}\subset\overline{{\cal G}_{(a,b)}}.$\\
B) If $c>a+b-1$ then ${\cal G}_{(c)}\cap\overline{{\cal G}_{(a,b)}}=\emptyset.$
\end{conj}

A result of the first author (\cite{E} Theorem~3) is that part A holds when $a-1$ divides $b-1$.
In~\cite{F4}, Furter proved both parts in the $a=2$ and $a=3$ case.   The method of \cite{F4} involves showing special cases of a new conjecture, the Rigidity Conjecture, which implies Conjecture~\ref{conj:dim2}.  Subsequently,  van den Essen and the first author \cite{EvdE} discovered a link between Furter's Rigidity Conjecture and the Factorial Conjecture, which is in turn related to the famous Jacobian Conjecture.

Of particular interest to us is the case $l=1$ and $m=3$. In~\cite{EF}, Furter and the first author proved that
$\overline{{\cal G}_{(11,3,3)}}$ intersects ${\cal G}_{(19)}$.  However, $\overline{{\cal G}_{(11,3,3)}}$ can not
contain ${\cal G}_{(19)}$ for dimensional reasons (as $19>11+3+3-1$). 
This implies that $\overline{{\cal G}_{(11,3,3)}}$ is not a union of some ${\cal G}_e$. More generally, we have (see~\cite{EF} Theorem~2):

\begin{theo}[Edo-Furter]\label{thm:sfpa2}
Let $a,b\ge 2$ and $c\ge 1$ be integers, and set $d=ab-1$.  Then
$${\cal G}_{(cd+a)}\cap\overline{{\cal G}_{(a+(c-1)d,b,a)}} \ne\emptyset$$
Moreover, if $(a,b)\ne(2,2)$ then
$${\cal G}_{(cd+a)}\not\subset\overline{{\cal G}_{(a+(c-1)d,b,a)}}.$$
\end{theo}

The reason why ${\cal G}_{(cd+a)}$ can not be a subset of $\overline{{\cal G}_{((c-1)d+a,b,a)}}$
in the case $(a,b)\ne(2,2)$ is the aforementioned dimension constraint.
Note that since $(c-1)d+a\le cd-1$, Theorem~\ref{thm:sfpa2} implies
$${\cal G}_{(cd+a)}\cap\overline{{\cal G}_{(cd-1,b,a)}} \ne\emptyset.$$
Since $\dim{\cal G}_{(cd+a)}\le\dim{\cal G}_{(cd-1,b,a)}-1$ (with equality if and only if $b=2$),
both assumptions of Question~\ref{ques:dim} are fulfilled and one can ask if the inclusion holds in this case. 
This is precisely the main result of this paper:

\renametheorem{Main Theorem}
\begin{namedtheorem}
Let $a,b\ge 2$ and $c\ge 1$ be integers, and set $d=ab-1$.  Then
$${\cal G}_{(cd+a)}\subset \overline{{\cal G}_{(cd-1,b,a)}}.$$
\end{namedtheorem}

To put this result in context, we make three observations.
\begin{enumerate}
\item This shows that the inclusion holds in the case  $(a,b)=(2,2)$ of Theorem~\ref{thm:sfpa2}.  That is, the necessary condition $(a,b) \neq (2,2)$ in the ``moreover'' statement of Theorem \ref{thm:sfpa2} is also sufficient.

\item When $b=2$, our main theorem is optimal in the dimensional sense; that is, 
$$\dim{\cal G}_{(a+c(2a-1))}=\dim{\cal G}_{(a,2,(2a-1)c-1)}-1$$
and gives ${\cal G}_{(a+(2a-1)c)}\subset \overline{{\cal G}_{(a,2,(2a-1)c-1)}}$.
With the notations of~\cite{F3}, we have: $(a+(2a-1)c)\sqsubseteq(a,2,(2a-1)c-1)$.
In~\cite{F2} and~\cite{F3}, Furter introduced an order denoted by $\preceq$ on ${\cal D}$ such that
$(d)\preceq(a,b)$ if and only if $d\le a+b-1$ and
$(d)\preceq(a,b,c)$ if and only if $d\le a+b+c-2$ (for all $a,b,c,d\ge 2$).
He conjectured (cf. Conjecture~7.1 in~\cite{F3}) that the orders $\sqsubseteq$ and $\preceq$ are equivalent.
Since $(a+(2a-1)c)\not\preceq(a,2,(2a-1)c-1)$, the main theorem gives a family of counterexamples to this conjecture.

\item Finally, as mentioned above, we know that
${\cal G}_{(a+b-1)}\subset\overline{{\cal G}_{(a,b)}}$ holds when $a-1$ divides $b-1$.
There seems to be a similar phenomenon for automorphisms with polydegree of length~$3$ with a $2$ in the middle: by the main theorem, the inclusion ${\cal G}_{(a+b+1)}\subset\overline{{\cal G}_{(a,2,b)}}$ holds when $2a-1$ divides $b+1$.  
\end{enumerate}

The remainder of this paper is devoted to proving the main theorem. The proof is quite technical
and is arranged in three parts: in section~2, we introduce the notion of {\em triangular polynomials}.
This new tool is central in the proof our result. In section~3, we compute a formal
inverse which arises naturally in the proof, and we prove that this formal inverse
is a triangular polynomial. In section~4, we use the results of the two previous sections to prove
the main theorem.  The reader may wish to begin with a light reading of section 4 in order to understand the motivation for the results of sections 2 and 3.

%%%%%%%%%%%%%%%%%%%%%%%
\section{Triangular polynomials}

\hspace*{.6cm} We fix a positive integer $a$ and $a+1$ variables $u_0,\ldots,u_a$.
Let $R$ denote the Laurent polynomial ring $\C[u_0,\ldots,u_{a-1}][u_a,u_a^{-1}]$, and set $U(Y)=\sum_{j=0}^au_jY^j$ in $R[Y]$.
Given an element $r\in R$ and $x=(x_0,\ldots,x_a)\in\C^a\times\C^*$ we denote by $r(x)\in\C$
the image of $r$ in the quotient $R/(u_0-x_0,\ldots,u_a-x_a)$,
which is canonically isomorphic to $\C$. We denote by $\Q_+$ the semiring of positive rational numbers.

\begin{defi}[Triangular polynomial]
Let $m\ge 1$ be an integer. We say that a polynomial $P(Y)\in R[Y]$ is $(m,U)$-\textit{triangular}
or simply $m$-\textit{triangular}
if $d:={\rm deg}(P(Y))\ge a$ and, writing $P(Y)=\sum_{l=0}^dp_lY^l$, we have for $d-a \leq l \leq d$,
$$p_l = q_l u_a^{m-1}u_{a-d+l}+P_l(u_{a-d+l+1},\ldots,u_a)$$
for some $q_l \in \Q_+$ and $P_l(u_{a-d+l+1},\ldots,u_a) \in \C[u_{a-d+l+1},\ldots,u_a]$.  Note that when $l=d$, this means $p_d=q_d u_a^m$ for some $q_d \in \Q_+$.
\end{defi}

\begin{ex} The polynomial $U(Y)$ is $1$-triangular.
\end{ex}

Triangular polynomials gain their utility from two useful properties.  First, in Proposition \ref{prop:surj}, we show that we can use them to generate a polynomial with the top $a+1$ coefficients specified.  Then, in the remainder of this section, we show that the set of triangular polynomials has some nice closure properties enabling us to build up other triangular polynomials from $U(Y)$.

\begin{prop}\label{prop:surj}
Let $t_{d-a},\ldots,t_d\in\C^*$ be nonzero complex numbers.
If the polynomial $P(Y)=\sum_{l=0}^dp_lY^l$ is $m$-triangular for some $m\ge 1$ of degree~$d$, then the map
$x\mapsto (t_{d-a}p_{d-a}(x),\ldots,t_dp_d(x))$ from $\C^a\times\C^*$ to
itself is surjective.
\end{prop}
\pf Composing by the bijection $(x_{d-a},\ldots,x_d)\mapsto(t_{d-a}x_{d-a},\ldots,t_dx_d)$
from $\C^a\times\C^*$ to itself, it's clear that we can assume that $t_i=1$ for all $i\in\{d-a,\ldots,d\}$.
Since $P(Y)$ is $m$-triangular there exist
$q_{d-a},\ldots,q_d\in\Q_+$ and $P_l(u_{a-d+l+1},\ldots,u_a)\in\C[u_{a-d+l+1},\ldots,u_a]$ for $l\in\{d-a,\ldots, d-1\}$
such that $p_d=q_du_a^m$ and $p_l=q_lu_a^{m-1}u_{a-d+l}+P_l(u_{a-d+l+1},\ldots,u_a)$
for all $l\in\{d-a,\ldots,d-1\}$. Given $y=(y_{d-a},\ldots,y_d)\in\C^a\times\C^*$,  we first choose
an $x_a \in \C^*$ such that $x_a^m=q_d^{-1}y_d$ (using the fact that $\C$ is algebraically closed).  To find $x_0,\ldots,x_{a-1}$, we note that the polynomial map $\theta = (p_{d-a},\ldots,p_{d-1})$ is a triangular automorphism of $\Spec \C[u_a,u_a^{-1}][u_0, \ldots, u_{a-1}]$ over $\Spec \C[u_a,u_a^{-1}]$.  Specializing to $u_a=x_a$, we obtain $\theta ^\prime$, a triangular automorphism of $\C^a = \Spec \C[u_0,\ldots,u_{a-1}]$.  Since $\theta ^\prime$ is a bijection of $\C^a$, we may set $x_i = (\theta ^\prime) ^{-1} (u_i)$ for $0 \leq i \leq a-1$.
This gives an antecedent $x=(x_0,\ldots,x_{a})$ to $y$. \qed

\begin{prop}\label{prop:tri}
Let $P(Y)\in R[Y]$ (resp. $Q(Y)\in R[Y]$) be an $m$-triangular (resp. $n$-triangular)
polynomial of degree $d$ (resp. $e$). Then:\\
1) If $d=e$ and $m=n$ then $P(Y)+Q(Y)$ is $m$-triangular.\\
2) $P(Y)Q(Y)$ is $m+n$-triangular.
\end{prop}
\pf The point  1) is clear since $\Q_+$ is closed under addition. Let us prove 2).
We write $P(Y)=\sum_{i=0}^dp_iY^i$ and $Q(Y)=\sum_{j=0}^eq_jY^j$ for some $p_i, q_j \in R$.
Then $P(Y)Q(Y)=\sum_{k=0}^{d+e}r_kY^k$ where $r_k=\sum_{i+j=k}p_iq_j \in R$.
Fix $k\in\{d+e-a,\ldots,d+e\}$, and let $i\in\{0,\ldots,d\}$ and $j\in\{0,\ldots,e\}$ be such that $i+j=k$.
Note that this implies $i\geq d-a$ and $j \geq e-a$.
To obtain  $r_k \in \Q^+u_a ^{m+n-1}u_{a-(d+e)+k}+\C[u_{a-(d+e)+k+1},\ldots,u_a]$, we will show that
\begin{enumerate}
\item $p_i q_j \in \C[u_{a-(d+e)+k+1},\ldots,u_a]$ when $i<d$ and $j<e$
\item $p_dq_{k-d}, p_{k-e}q_e \in \Q_+ u_a ^{m+n-1}u_{a-(d+e)+k} +\C[u_{a-(d+e)+k+1},\ldots,u_a] $
\end{enumerate}

So we first assume $i<d$ and $j < e$.  Then since $P(Y)$, $Q(Y)$ are triangular, we have
\begin{align*}
p_i &\in\Q_+ u_a^{m-1}u_{a-d+i}+\C[u_{a-d+i+1},\ldots,u_a]\\
q_j &\in\Q_+ u_a^{n-1}u_{a-e+j}+\C[u_{a-e+j+1},\ldots,u_a]
\end{align*}
and since $i=k-j>k-e$ and $j=k-i>k-d$, we see $p_i q_j \in \C[u_{a-(d+e)+k+1},\ldots,u_a]$.

If $i=d$  we have $p_d \in \Q_+ u_a ^m$, and thus $p_d q_{k-d} \in \Q_+ u_a ^{m+n-1}u_{a-(d+e)+k}+\C[u_{a-(d+e)+k+1},\ldots,u_a]$.  The $j=e$ case is similar.  \qed\\

\begin{ex}
Let $b\ge 1$ be an integer. The polynomial $U(Y)^b$ is $b$-triangular.
\end{ex}

It will be convenient to have the following definition of {\em linear polynomials}, which behave like $1$-triangular polynomials.
\begin{defi}[linear polynomial]
We say that a polynomial $P(Y)\in R[Y]$ is $U$-\textit{linear} or simply \textit{linear}
if $0\le d:={\rm deg}(P(Y))\le a$ and if $P(Y)=\sum_{l=0}^dp_lY^l$ with
$p_l\in\Q_+u_{a-d+l}$ for all $l\in\{0,\ldots,d\}$.
\end{defi}
\begin{ex}
Obviously, if $P(Y)\in\Q_+$ then $P(Y)$ is linear.
\end{ex}
\begin{ex} $U(Y)$ and $U^\prime(Y)$ are $U$-linear.
\end{ex}

\begin{prop}\label{prop:lin}
Let $P(Y)\in R[Y]$ be $m$-triangular and let $Q(Y) \in R[Y]$ be a linear polynomial.\\
0) For all $q\in\Q_+$, $qP(Y)$ is $m$-triangular, and $qQ(Y)$ is linear.\\
1) If ${\rm deg}(P(Y))\ge a+1$ then $P'(Y)$ is $m$-triangular.\\
2) If $Q(Y)$ is a nonconstant polynomial, then $Q'(Y)$ is linear.\\
3) $P(Y)Q(Y)$ is $m+1$-triangular.
\end{prop}
\pf
We write $P(Y)=\sum_{i=0}^dp_iY^i$ and $Q(Y)=\sum_{j=0}^eq_jY^j$ where $d:={\rm deg}(P(Y))$ and $e:={\rm deg}(Q(Y))$.  The assertion 0)  is obvious.  For point 1), we assume that ${\rm deg}(P(Y))\ge a+1$.  Then, we have
$$P'(Y)=\sum_{i=0}^{d-1}(i+1) p_{i+1} Y^i .$$
Since $P(Y)$ is $m$-triangular, we have $$(i+1)p_{i+1} \in \Q_+u_a^{m-1}u_{a-d+i+1}+\C[u_{a-d+i+1+1},\ldots,u_a]$$
for $i\in\{d-1-a,\ldots,d-1\}$.
Hence $P'(Y)$ is  $m$-triangular.

To prove 2), suppose $\deg Q(Y) = e>0$.  Then $$Q'(Y)=\sum_{j=0}^{e-1}(j+1) q_{j+1} Y^j.$$  But since $q_{j} \in \Q_+ u_{a-d+j}$ for each $0 \leq j \leq e$, we have $(j+1)q_{j+1} \in \Q_+ u_{a-(d-1)+j}$ for each $ 0 \leq j \leq e-1$.  Thus $Q^\prime(Y)$ is linear.

Finally, we remark that the proof of part 3) is similar to the proof of point 2) of the preceding proposition. \qed \\

Using Proposition~\ref{prop:lin} inductively we deduce:
\begin{coro}\label{coro:lin}
Let $P(Y)\in R[Y]$ be an $m$-triangular polynomial of degree~$d$.
Let $k_0,\ldots,k_a\ge 0$ be integers. Then the polynomial
$$P(Y)\prod_{j=0}^a(U^{(j)}(Y))^{k_j}$$
is $m+\sum_{j=0}^ak_j$-triangular. In particular, if $n=\sum_{j=0}^ak_j$ and
$\sum_{j=0}^ajk_j=n$, then this polynomial is $m+n$-triangular and has degree $d+(a-1)n$.
\end{coro}

%%%%%%%%%%%%%%%%%%%%%%%
\section{A formal inverse computation}

Let $a,b\ge 2$ be integers. Let $R$ be a $\Q$-algebra, and let $u_0,\ldots,u_a$ be in $R$.
We consider the polynomial $U(Y)=\sum_{j=0}^au_jY^j\in R[Y]$.
We denote by $I(Y,Z)$ the formal inverse of $Y+ZU(Y)^b$ in $R[[Z]][Y]$.
In other words, $I(Y,Z)$ is the unique element of $R[[Z]][Y]$ satisfying $I(Y+ZU(Y)^b,Z)=Y$.
We denote by $v_k(Y) \in R[Y]$ the coefficient of $Z^k$ in $U(I(Y,Z))$ for all $k\ge 0$, so that
$$U(I(Y,Z))=\sum_{k=0}^{\infty}v_k(Y)Z^k.$$
Plugging in $Y+ZU(Y)^b$ for $Y$, we obtain
\begin{align*}
U(Y) & =U(I(Y+ZU(Y)^b,Z)) =\sum_{k=0}^{\infty}v_k(Y+ZU(Y)^b)Z^k .
\end{align*}
Now, applying Taylor's formula and setting $m=j+k$, we deduce
\begin{align*}
U(Y) &=\sum_{k=0}^{\infty}Z^k\sum_{j=0}^{\infty}{v_k^{(j)}(Y)\over j!}\,Z^jU(Y)^{bj} =\sum_{m=0}^{\infty}Z^m\sum_{j=0}^m{U(Y)^{bj}\over j!}\,v_{m-j} ^{(j)}(Y).
\end{align*}
Since $U(Y)$ has $Z$-degree zero, we obtain the following recursive relation for $v_k(Y)$:
\begin{align*}
v_0(Y) & =U(Y) \\
v_m(Y) &= - \sum _{j=1}^m\frac{U(Y)^{bj}}{j!} v_{m-j} ^{(j)}(Y) &  m\ge 1
\end{align*}

We would like to find a non-recursive formula for $v_k(Y)$ (see Theorem \ref{thm:fi}).  To this end, we define for any integer $\lambda \geq 0$:
\begin{align}
w_{0,\lambda}(Y) &= 1/(\lambda+1) \notag \\
w_{n,\lambda}(Y) &= (\lambda-n+2)U'(Y)w_{n-1,\lambda}(Y)+U(Y)w_{n-1,\lambda}'(Y) & n\ge 1 \label{wnlambda}
\end{align}

These polynomials arise naturally in our computation of $v_k(Y)$.  However, as they are also defined recursively, we first find a non-recursive formula for $w_{n, \lambda}(Y)$.
\begin{prop}\label{prop:w}
Let $1\le n\le\lambda+1$ be an integer. We denote by $I_n$ the set of sequences $(k_0,\ldots,k_a)$ of non-negative integers such that
$\sum_{j=0}^ak_j=n$ and $\sum_{j=0}^ajk_j=n$. Then
$$w_{n,\lambda}(Y)=\sum_{(k_0,\ldots,k_a)\in I_n}q_{(k_0,\ldots,k_a)}\prod_{j=0}^a(U^{(j)}(Y))^{k_j}$$
for some $q_{(k_0,\ldots,k_a)}\in \Q_+$.
\end{prop}
\pf We induct on $n$.  The result is trivial when $n=0$. Let $n\ge 1$ be an integer, and assume the formula holds for $w_{n-1,\lambda}(Y)$.
For all integers $0\le j\le a$, we denote by $e_j$ the $j$-th standard basis vector of $\Q^{a+1}$.
If $(k_0,\ldots,k_a)$ belongs to $I_n$, then
\begin{equation}\label{w1}
U'(Y)\prod_{j=0}^a(U^{(j)}(Y))^{k_j}=\prod_{j=0}^a(U^{(j)}(Y))^{l_j}
\end{equation}
where $(l_0,\ldots,l_a)=(k_0,\ldots,k_a)+e_1\in I_{n+1}$ and
\begin{equation}\label{w2}
U(Y)\left(\prod_{j=0}^a(U^{(j)}(Y))^{k_j}\right)'=\sum_{j=0}^{a-1}k_j\prod_{i=0}^a(U^{(j)}(Y))^{l_{i,j}}
\end{equation}
where $(l_{0,j},\ldots,l_{a,j})=(k_0,\ldots,k_a)+e_0-e_j+e_{j+1}\in I_{n+1}$, for all integers $0\le j\le a-1$.
Combining \eqref{w1}, \eqref{w2}, \eqref{wnlambda}, and the induction hypothesis, we see
\begin{align*}
w_{n,\lambda} &=  \sum _{(k_0,\ldots,k_a) \in I_{n-1}} q_{k_0,\ldots,k_a} \left( (\lambda -n +2) U^\prime(Y)  \prod _{j=0} ^a (U^{(j)}(Y))^{k_j} + U(Y) \left( \prod _{j=0} ^a (U^{(j)}(Y))^{k_j}\right) ^\prime \right)\\
&=\sum _{(k_0,\ldots,k_a) \in I_{n-1}}q_{k_0,\ldots,k_a} \left( (\lambda -n +2)  \prod _{j=0} ^a (U^{(j)}(Y))^{l_{0,j}} + \sum _{j=0} ^{a-1} k_j \prod _{i=0} ^a (U^{(j)}(Y)) ^{l_{i,j}} \right)\\
&= \sum _{(k_0,\ldots,k_a) \in I_n} \tilde{q}_{(k_0,...,k_a)} \prod _{j=0} ^a (U^{(j)}(Y))^{k_j}
\end{align*}
Since $\lambda > n-2$ and $q_{k_0,\ldots,k_a}, k_j \in \Q_+$, we have $\tilde{q}_{(k_0,\ldots,k_a)} \in \Q_+$, completing the proof. \qed

\begin{lem}\label{lem:fi}
Let $n,r$ be non-negative integers, and let $k,m$ be positive integers with $1 \leq k \leq m$. Then $$\sum _{j=0} ^m (-1)^j {m \choose j} w_{k, (m+r-j)b} ^{(n)}(Y) =0.$$
\end{lem}
\pf We denote by $S(n,k,m,r)$ the sum in the lemma. It's clear that given integers $1 \leq k \leq m$ and $r \geq 0$,
$S(0,k,m,r)=0$ implies $S(n,k,m,r)=0$ for all $n\ge 0$ (differentiate $n$ times).
We thus assume $n=0$ and prove $S(0,k,m,r)=0$ for all $m\geq k$ and $r \geq 0$ by induction on $k$.

First, if $k=1$, we have, for any $m\geq k$, $r \geq 0$ (noting that $U_{1,\lambda}(Y)=U^\prime(Y)$ for any $\lambda$),
$$\sum _{j=0} ^m (-1) ^j {m \choose j} w_{1,\lambda} (Y) = {U^\prime(Y)} \sum _{j=0} ^{m} (-1)^j {m \choose j}= U^\prime(Y) (1-1)^m=0.$$

Now, suppose $k > 1$. Let $m \geq k$, $r \geq 0$ be integers, and assume $S(0,k-1,m',r')=0$ for all $m' \geq k$, $r' \geq 0$.
We set $\lambda_j =(m+r-j)b$ and we observe
\begin{align}
\lambda_j-k+2 &=b(m-j)+(rb-k+2) \label {lambdajk2} \\
\lambda _j &=((m-1)+(r+1)-j)b . \label{lambdaj}
\end{align}
We now compute by \eqref{wnlambda} and \eqref{lambdajk2}
\begin{align*}
S(0,k,m,r) &= \sum _{j=0} ^m (-1) ^j {m \choose j} \left( \left(\lambda_j-k+2\right)U'(Y)w_{k-1,\lambda _j}(Y)+U(Y)w_{k-1,\lambda _j}'(Y) \right)\\
           &= bU'(Y) \sum _{j=0} ^m (-1) ^j {m \choose j} (m-j) w_{k-1,\lambda _j}(Y) + \\
 &\phantom{xxx} (rb-k+2)U'(Y)\sum _{j=0} ^m (-1)^j {m \choose j}w_{k-1,\lambda _j}(Y)+\\
&\phantom{xxx} U(Y) \sum _{j=0} ^m (-1) ^j {m \choose j} w_{k-1,\lambda _j}'(Y) \\
\end{align*}
Noting that ${m \choose j}(m-j) = m{m-1 \choose j}$, we see from \eqref{lambdaj} that the first sum is $mbU^\prime(y)S(0,k-1,m-1,r+1)$.  Thus, we have
\begin{align*}
S(0,k,m,r)&=mbU'(Y)S(0,k-1,m-1,r+1) +(rb-k+2)U'(Y)S(0,k-1,m,r)\\
&\phantom{xxx}+U(Y)S(1,k-1,m,r)
\end{align*} 
By the induction hypothesis, each of these three terms is zero, so $S(0,k,m,r)=0$.  \qed

\begin{theo}\label{thm:fi}
For $m \geq 0$, $n \geq 0$, we have
$$v_m ^{(n)}(Y) = \frac{(-1)^m}{m!} U(Y)^{bm-m-n+1}\, w_{m+n,bm}(Y).$$
\end{theo}
\pf The proof is by induction on $(m,n)$.  First, we  verify the $(0,0)$ case:
$$v_0(Y)=U(Y)=\frac{(-1)^0}{0!}U(Y)w_{0,0}(Y)$$
Let $m \geq 0$, $n \geq 1$ be integers. We assume the $(m,n-1)$ case and we prove the $(m,n)$ case.
We set $N=m+n$ and $\lambda=bm$ and we use Leibniz's rule.
\begin{align*}
v_{m}^{(n)}(Y) &=\frac{\partial}{\partial Y} v_m ^{(n-1)}(Y)
=\frac{\partial}{\partial Y} \left(\frac{(-1)^m}{m!} U(Y)^{\lambda-N+2} w_{N-1,\lambda}(Y)\right) \\
&\hspace*{-1cm}=\frac{(-1)^m}{m!} \left( (\lambda-N+2)U(Y)^{\lambda-N+1}U'(Y)w_{N-1,\lambda}(Y)+U(Y)^{\lambda-N+2}w_{N-1,\lambda}'(Y)\right)\\
&\hspace*{-1cm}=\frac{(-1)^m}{m!} U(Y)^{\lambda-N+1} \left((\lambda-N+2)U'(Y)w_{N-1,\lambda}(Y)+U(Y)w_{N-1,\lambda}'(Y)\right) \\
&\hspace*{-1cm}=\frac{(-1)^m}{m!} U(Y)^{\lambda-N+1} w_{N,\lambda}(Y)
\end{align*}

Let $m \geq 1$ be an integer. We assume the $(m',n)$ case for all pairs $(m',n)$ such that $m'<m$
and we prove the $(m,0)$ case.
\begin{align*}
v_m (Y)&= -\sum _{j=1} ^m \frac{U(Y)^{bj}}{j!} v_{m-j} ^{(j)} (Y)\\
&=-\sum _{j=1} ^m \frac{U(Y)^{bj}}{j!} \left( \frac{(-1)^{m-j}}{(m-j)!} U(Y)^{b(m-j)-(m-j)-j+1} w_{m,b(m-j)}(Y)\right) \\
&=-\frac{(-1)^m}{m!} U(Y)^{bm-m+1} \sum _{j=1} ^m (-1)^{j}{m \choose j} w_{m,b(m-j)} (Y)\\
&= \frac{(-1)^m}{m!} U(Y)^{bm-m+1} w_{m,bm}(Y).
\end{align*}
The last equality following from Lemma~\ref{lem:fi} ($S(0,m,m,0)=0$). \qed

\begin{theo}\label{thm:triangular}
Let $m\ge 0$ be an integer. The polynomial $(-1)^m\,m!\,v_m(Y)$ is $(U,bm+1)$-triangular
of degree $(ab-1)m+a$.
\end{theo}
\pf By Theorem~\ref{thm:fi}, we have $(-1)^m\,m!\,v_m(Y)=U(Y)^{bm-m+1}\, w_{m,bm}(Y).$
Using Proposition~\ref{prop:w}, we see that this polynomial is a sum of terms (not all equal to zero)
of the form $U(Y)^{bm-m+1}\prod_{j=0}^a(U^{(j)}(Y))^{k_j}$ where $k_0,\ldots,k_a\ge 0$ are integers
such that $\sum_{j=0}^ak_j=m$ and $\sum_{j=0}^ajk_j=m$. By 2) of Proposition~\ref{prop:tri}, $P(Y)=U(Y)^{bm-m+1}$
is $bm-m+1$-triangular of degree $a(bm-m+1)$. By Corollary \ref{coro:lin}, we deduce that
all terms $U(Y)^{bm-m+1}\prod_{j=0}^a(U^{(j)}(Y))^{k_j}$ are $bm+1$-triangular
of degree $(ab-1)m+a$ and we conclude that the sum is $bm+1$-triangular by 1) of Proposition~\ref{prop:tri}. \qed

%%%%%%%%%%%%%%%%%%%%%%%
\section{Main theorem}

In this section, we prove:
\begin{namedtheorem}
Let $a,b\ge 2$ and $c\ge 1$ be integers, we set $d=ab-1$, we have:
$${\cal G}_{(cd+a)}\subset \overline{{\cal G}_{(cd-1,b,a)}}.$$
\end{namedtheorem}

\pf Let $\tau\in {\cal G}_{(cd+a)}$ be a triangular automorphism of degree $cd+a$. Then we can write
$$\tau=(rX+\sum_{j=0}^{cd+a}y_jY^j,sY+t)$$
for some $r,y_{cd+a},s\in\C^*$ and $y_0,\ldots,y_{cd+a-1},t\in\C$. To prove $\tau\in\overline{{\cal G}_{(cd-1,b,a)}}$, using
the Valuation Criterion (see Theorem~\ref{thm:valcri}), we construct an automorphism $\s_Z\in {\cal G}_{(cd-1,b,a)}(\C[Z])$
such that writing $\s_Z=(f_Z,g_Z)$ , then going modulo $Z$, $f_Z\equiv rX+\sum_{j=0}^{cd+a}y_jY^j$ and $g_Z\equiv sY+t$.

We continue to use the notations of section~2 and~3. By Theorem~\ref{thm:triangular},
$(-1)^cc!v_c(Y)$ is $(U,bm+1)$-triangular of degree $cd+a$.
We write $v_c(Y)=\sum_{l=0}^{cd+a}p_lY^l$ with $p_l\in R$ for $l\in\{0,\ldots,cd+a\}$.
We apply Proposition~\ref{prop:surj}, with $t_l={r\over (-1)^cc!}$ for $l\in\{cd,\ldots,cd+a\}$
and $P(Y)=(-1)^cc!v_c(Y)$: there exists $x\in\C^a\times\C^*$, such that $rp_l(x)=y_l$ for $l\in\{cd,\ldots,cd+a\}$.

We now fix such an $x$.  Given a polynomial $Q(Y,Z)=\sum_{l,m}q_{l,m}Y^lZ^m$ in $R[Y,Z]$ we denote by
$\overline{Q}(Y,Z)$ the specialization $\sum_{l=0}^eq_{l,m}(x)Y^lZ^m\in\C[Y,Z]$. The map $Q(Y,Z)\mapsto\overline{Q}(Y,Z)$ is a ring homomorphism
from $R[Y,Z]$ to $\C[X,Y]$ and if $Q(Y)\in R[Y]$ then $\deg(\overline{Q}(Y))=\deg(Q(Y))$. We have: $\overline{U}(Y)=\sum_{j=0}^ax_jY^j$
and $r\overline{v_c}(Y)=\sum_{j=0}^{cd+a}y_jY^j+E(Y)$ where $E(Y)\in\C[Y]$ is such that $\deg(E(Y))\le cd-1$.

We truncate the power series $\overline{U}(\overline{I}(Y,Z))=\sum_{k=0}^{\infty}\overline{v_k}(Y)Z^k$ to the polynomial  $V(Y,Z)=\sum_{k=0}^{c-1}\overline{v_k}(Y)Z^k$.
We consider the following three triangular automorphisms of $\C(Z)[X,Y]$:
\begin{align*}
\tau_1 & =(Z^cX+\overline{U}(Y),Y) \\ \tau_2&=(X+ZY^b,Y) \\ \tau_3&=(rZ^{-c}(X-V(Y,Z))-E(Y)+ZY^{cd-1},sY+t)
\end{align*}

We set $\s_Z=\tau_3\pi\tau_2\pi\tau_1$ where $\pi=(Y,X)$. We then have, letting $W(Y,Z)=Y+Z(Z^cX+\overline{U}(Y))^b$ :
\begin{align*}
\s_Z(X)&=r\left(X+Z^{-c}\left[\,\overline{U}(Y)-V(W(Y,Z),Z)\,\right]-E(W(Y,Z))+ZW(Y,Z)^{cd-1}\right)\\
\s_Z(Y)&=sW(Y,Z)+t
\end{align*}

Focusing on $\s_Z(X)$, we compute modulo $Z^{c+1}$:
\begin{align*}
\overline{U}(Y)-V(W(Y,Z),Z) &\equiv\overline{U}(Y)-\overline{U}(\overline{I}(W(Y,Z),Z))+\overline{v_c}(W(Y,Z))Z^c  \\
                 &\equiv\overline{U}(Y)-\overline{U}(\overline{I}(Y+Z\overline{U}(Y)^b,Z))+\overline{v_c}(Y)Z^c  \\
 &\equiv\overline{v_c}(Y)Z^c
\end{align*}
This computation proves that $\s_Z(X)\in\C[Z][X,Y]$.  Since the Jacobian of $\s_Z$ is $1$,
we deduce from the overring principle (see Lemma~1.1.8 p. 5 in~\cite{vdE}) that $\s_Z$ is a $\C[Z]$-automorphism of $\C[Z][X,Y]$.
Moreover, we have $$\lim_{Z\to 0}\s_Z = (r(X+\overline{v_c}(Y))-E(Y),sY+t)=\tau.$$
It's clear that $\deg(\tau_1)=\deg(\overline{U}(Y))=a$ and $\deg(\tau_2)=b$.
Since $\deg(E(Y))\le cd-1$ and for all $k\le c-1$, $\deg \overline{v_k}(Y)=kd+a\le(c-1)d+a\le cd-1$,
we have $\deg(\tau_3)=cd-1$ and thus $\s_Z\in{\cal G}_{(cd-1,b,a)}(\C[Z])$. \qed

\end{document}